\documentclass[12pt]{article}
\usepackage{epic,eepic,epsf,epsfig}
\usepackage{amsmath,enumerate}
\usepackage{amssymb}
\usepackage{amsbsy}

\usepackage{amsfonts}

 \def\Omega{\Omega}

\def\f{\noindent}

\def\Syl{\hbox{\rm Syl}}

\newcommand{\qed}{\mbox{\raisebox{0.7ex}{\fbox{}}} \vspace{4truemm}}
\def\demo{\f {\bf Proof.}\hskip10pt}

\begin{document}

\baselineskip 16pt

\title{ \vspace{-1.2cm}
Finite groups in which every maximal subgroup is nilpotent or normal or has $p'$-order
\thanks{\scriptsize
J. Shi was supported by Shandong Provincial Natural Science Foundation, China (ZR2017MA022 and ZR2020MA044) and NSFC (11761079). R. Shen 
was supported by NSFC (12161035).
\newline
 \hspace*{0.5cm} \scriptsize $^{\ast\ast}$ Corresponding author.
\newline
\hspace*{0.5cm} \scriptsize{E-mail addresses: jiangtaoshi@126.com\,(J. Shi),\,\,ln18865550588@163.com\,(N. Li),\,\,shenrulin@hotmail.com\,(R. Shen).}}}

\author{Jiangtao Shi$^{a,\,\ast\ast}$,\,\,Na Li$^b$,\,\,Rulin Shen$^c$\\
\\
{\small{\em $^a$ School of Mathematics and Information Sciences, Yantai University, Yantai 264005, P.R. China}}\\
{\small{\em $^b$ Department of Mathematics, Beijing Jiaotong University, Beijing 100044, P.R. China}}\\
{\small{\em $^c$ Department of Mathematics, Hubei Minzu University, Enshi 445000, P.R. China}}}
\date{ }

\maketitle \vspace{-.8cm}

\begin{abstract}
Let $G$ be a finite group and $p$ a fixed prime divisor of $|G|$. Combining the nilpotence, the normality and the order of
groups together, we prove that if every maximal subgroup of $G$ is nilpotent or normal or has $p'$-order, then
(1) $G$ is solvable; (2) $G$ has a Sylow tower; (3) There exists at most one prime divisor $q$ of $|G|$ such that
$G$ is neither $q$-nilpotent nor $q$-closed, where $q\neq p$.

\medskip \f {\bf Keywords:} maximal subgroup; nilpotent; normal; $p'$-order; $q$-nilpotent; $q$-closed
\end{abstract}

\section{Introduction}

In this paper all groups are assumed to be finite. It is known that a group $G$ is nilpotent or minimal non-nilpotent
if every maximal subgroup of $G$ is nilpotent, and a group $G$ is nilpotent if every maximal subgroup of $G$ is normal,
see {\rm\cite[Theorem 9.1.9]{rob}} and {\rm\cite[Theorem 5.2.4]{rob}}, respectively. As a generalization, combining
the nilpotence and the normality of groups together, Li and Guo {\rm\cite[Theorem 1.2]{qian}} proved that if all
non-normal maximal subgroups of a group $G$ are nilpotent then $G$ is solvable and $G$ is $p$-nilpotent for some prime
$p$, that is, $G$ has a normal $p$-complement. It is clear that the hypothesis that all non-normal maximal subgroups of
a group $G$ are nilpotent is equivalent to the hypothesis that every maximal subgroup of $G$ is nilpotent or normal.
Lu, Pang and Zhong {\rm\cite[Theorems 1.3 and 3.5]{lu}} proved that if every maximal subgroup of a group $G$ is nilpotent
or normal then $G$ is solvable and $G$ is $p$-nilpotent and $q$-closed for some primes $p$ and $q$, that is, $G$ has a
normal $p$-complement and the Sylow $q$-subgroup of $G$ is normal. In {\rm\cite[Theorem 1.1]{li}} we gave an elementary
proof of the solvability of a group in which every maximal subgroup is nilpotent or normal. Moreover, the first
author of this paper {\rm\cite[Theorem 5]{shi}} proved that such a group has a Sylow tower.

The order of subgroups or groups play an important role in characterizing the solvability of groups. Feit-Thompson
theorem shows that a group having odd order is solvable. Thompson {\rm\cite[Theorem 10.4.2]{rob}} implies that a group
having a nilpotent maximal subgroup of odd order is also solvable. Moreover, it is easy to see that a group $G$ satisfying
$(|G|,15)=1$ is solvable by the classification of minimal non-abelian simple groups.

Note that the nilpotence, the normality and the order are three distinct characteristic properties of groups.
In this paper, we combine the nilpotence, the normality and the order of groups together to give a complete
characterization of the structure of the group in which every maximal subgroup is nilpotent or normal or has
$p'$-order for a fixed prime divisor $p$ of its order, which extends and generalizes the researches
in {\rm\cite[Theorem 1.2]{qian}}, {\rm\cite[Theorems 1.3 and 3.5]{lu} and  {\rm\cite[Theorem 5]{shi}}.

First we obtain a basic structural property of the non-solvable groups, see Theorem 1.1, whose proof is given in Section~\ref{th1}.

\medskip

\f {\bf Theorem 1.1}\ \ {\it Let $G$ be a non-solvable group and $p$ a fixed prime divisor of $|G|$, then $G$ has non-nilpotent
maximal subgroups of order divisible by $p$.}

\medskip

\f {\bf Remark 1.2}\ \ {\rm Note that 2 is a prime divisor of the order of the alternating group $A_5$ but $A_5$ has no
non-nilpotent maximal subgroups of 2'-order. This example shows that a non-solvable group might not have non-nilpotent
maximal subgroups of $p'$-order for some fixed prime divisor $p$ of its order.}

\medskip

Our main results are the following Theorems 1.3, 1.4 and 1.7, whose proofs are given in Sections~\ref{th3}, ~\ref{th4}
and ~\ref{th7}, respectively.

\medskip

\f {\bf Theorem 1.3}\ \ {\it Let $G$ be a group and $p$ a fixed prime divisor of $|G|$. If every maximal subgroup of $G$ is
nilpotent or normal or has $p'$-order, then $G$ is solvable.}

\medskip

\f {\bf Theorem 1.4}\ \ {\it Let $G$ be a group and $p$ a fixed prime divisor of $|G|$. If every maximal subgroup
of $G$ is nilpotent or normal or has $p'$-order, then $G$ has a Sylow tower.}

\medskip

\f {\bf Remark 1.5}\ \ {\rm Let $G$ be a group and $p$ a fixed prime divisor of $|G|$. If we assume that every maximal
subgroup of $G$ is nilpotent or normal or has order divisible by $p$, we cannot get that $G$ has a Sylow tower.
For example, every maximal subgroup of the alternating group $A_5$ has order divisible by 2 then $A_5$ naturally
satisfies the hypothesis. But $A_5$ is a non-solvable group which has no Sylow tower.}

\medskip

\f {\bf Remark 1.6}\ \ {\rm The alternating group $A_4$ implies that the Sylow tower of the group $G$ in Theorem
1.4 might not be supersolvable type.}

\medskip

\f {\bf Theorem 1.7}\ \ {\it Let $G$ be a group and $p$ a fixed prime divisor of $|G|$. If  every maximal subgroup of
$G$ is nilpotent or normal or has $p'$-order, then there exists at most one prime divisor $q$ of $|G|$ such that $G$
is neither $q$-nilpotent nor $q$-closed, where $q\neq p$.}

\medskip

Comparing with Theorem 1.7, we have the following three results, all of whose proofs are given in Section~\ref{th8}.

\medskip

\f {\bf Theorem 1.8}\ \ {\it Let $G$ be a group in which every maximal subgroup is nilpotent or normal, then
$G$ is either $q$-nilpotent or $q$-closed for each prime divisor $q$ of $|G|$.}

\medskip

\f {\bf Theorem 1.9}\ \ {\it Let $G$ be a group and $p$ a fixed prime divisor of $|G|$. If every maximal subgroup of $G$
is nilpotent or has $p'$-order, then $G$ is either $q$-nilpotent or $q$-closed for each prime divisor $q$ of $|G|$.}

\medskip

\f {\bf Theorem 1.10}\ \ {\it Let $G$ be a group and $p$ a fixed prime divisor of $|G|$. If every maximal subgroup of $G$ is
normal or has $p'$-order, then $G$ is either $q$-nilpotent or $q$-closed for each prime divisor $q$ of $|G|$.}

\section{Proof of Theorem 1.1}\label{th1}

\medskip

\demo Let $G$ be a counterexample of minimal order. Then every maximal subgroup of $G$ is nilpotent or has
$p'$-order. Note that $p$ is a prime divisor of $|G|$, one has that $G$ must have maximal subgroups of order
divisible by $p$ and such maximal subgroups are nilpotent by the hypothesis. Since $G$ is non-solvable and
both a nilpotent group and a minimal non-nilpotent group are solvable by {\rm\cite[Theorem 9.1.9]{rob}},
$G$ must have non-nilpotent maximal subgroups which have $p'$-order.

Let $P\in\Syl_p(G)$. First suppose that there exists a nontrivial subgroup $U$ of $P$ such that $U\unlhd G$.
Considering the quotient group $G/U$. (1) Assume $U<P$. Then $G/U$ is a non-solvable group of order divisible
by $p$ since $G$ is non-solvable and $U$ is solvable. By the minimality of $G$, one has that $G/U$ has a
non-nilpotent maximal subgroup $H/U$ of order divisible by $p$. It follows that $H$ is a non-nilpotent maximal
subgroup of $G$ of order divisible by $p$, a contradiction. (2) Assume $U=P$. Then $U$ is a Sylow $p$-subgroup
of $G$. Let $M$ be a non-nilpotent maximal subgroup of $G$ of $p'$-order. It is clear that $U\nleq M$ and
$U\cap M=1$. Then $G=U\rtimes M$, the semidirect product of $U$ and $M$. For every maximal subgroup $M_0$ of
$M$, one has that $U\rtimes M_0$ is a maximal subgroup of $G$ of order divisible by $p$. By the hypothesis,
$U\rtimes M_0$ must be nilpotent. It follows that every maximal subgroup of $M$ is nilpotent, which implies
that $M$ is a nilpotent group or a minimal non-nilpotent group. One has that $M$ is solvable. Then
$G=U\rtimes M$ is solvable, this contradicts that $G$ is non-solvable.

Next suppose that every nontrivial subgroup $U$ of $P$ is not normal in $G$, that is, $N_G(U)<G$. Since the
order of $N_G(U)$ is divisible by $p$, $N_G(U)$ can only be contained in a nilpotent maximal subgroup of $G$
of order divisible by $p$. Then $N_G(U)$ is nilpotent. By {\rm\cite[IV, Theorem 5.8(b)]{huppert}}, one has
that $G$ is $p$-nilpotent. There exists a normal subgroup $K$ of $G$ such that $G=K\rtimes P$. Let $L$ be
a maximal subgroup of $G$ of $p'$-order. Assume $K\nleq L$. Then $KL>L$. One has $G=KL$, which implies
that $G$ has $p'$-order, a contradiction. Thus $K\leq L$. Since $K$ is a $p'$-Hall subgroup of $G$ and
$L$ has $p'$-order, one has that $K=L$ is a maximal subgroup of $G$. It follows that $P\cong G/K$ is a
cyclic group of order $p$. Note that $P<G$ and $P$ can only be contained in a nilpotent maximal subgroup
of $G$. (1) Assume $p=2$. By {\rm\cite[IV, Theorem 7.4]{huppert}}, one has that $G$ is solvable, a
contradiction. (2) Assume $p>2$. $(i)$ Suppose $Z(G)\neq 1$, that is, the center of $G$ is not equal
to 1. Considering the quotient group $G/Z(G)$. If $p\nmid|G/Z(G)|$, then $P\leq Z(G)$. It follows that
$G=K\times P$. Note that $K\cong G/P$ is non-solvable, one has that $K$ has a non-nilpotent maximal
subgroup $K_0$. Then $K_0\times P$ is a non-nilpotent maximal subgroup of $G$ of order divisible by
$p$, a contradiction. If $p\mid|G/Z(G)|$, then $G/Z(G)$ is a non-solvable group of order divisible
by $p$. By the minimality of $G$, $G/Z(G)$ has a non-nilpotent maximal subgroup $R/Z(G)$ of order
divisible by $p$. It follows that $R$ is a non-nilpotent maximal subgroup of $G$ of order divisible
by $p$, a contradiction, too. $(ii)$ Suppose $Z(G)=1$. By $\rm\cite[Theorem 1]{rose}$, every nilpotent
maximal of $G$ is a Sylow 2-subgroup of $G$. It implies that $P$ cannot be contained in any nilpotent
maximal subgroup of $G$ since $p>2$, this contradicts that $P$ can only be contained in a nilpotent
maximal subgroup of $G$.

Thus the counterexample of minimal order does not exist and so $G$ has non-nilpotent maximal
subgroups of order divisible by $p$.\hfill\qed

\section{Proof of Theorem 1.3}\label{th3}

\medskip

\demo Let $G$ be a counterexample of minimal order. Since $G$ is non-solvable, $G$ has non-nilpotent
maximal subgroups of order divisible by $p$ by Theorem 1.1. Then by the hypothesis, the group $G$ is
in particular not simple.

Let $N$ be a minimal normal subgroup of $G$. We will show that $N$ is non-solvable. (1) Suppose that $G/N$
has $p'$-order. Then $N$ contains the Sylow $p$-subgroup of $G$. For every non-nilpotent maximal
subgroup $L/N$ of $G/N$, $L$ is a non-nilpotent maximal subgroup of $G$ of order divisible by $p$.
By the hypothesis, $L\unlhd G$. It follows that $G/N$ is a group in which every maximal subgroup is
nilpotent or normal. By {\rm\cite[Theorem 1.1]{li}}, $G/N$ is solvable. Thus $N$ is non-solvable since
$G$ is non-solvable. (2) Suppose that $G/N$ has order divisible by $p$. Since the hypothesis of the
theorem holds for $G/N$ and $|G/N|<|G|$, one has that $G/N$ is solvable by the minimality of $G$.
Thus we also have that $N$ is non-solvable.

It follows that $\it\Phi$$(G)=1$ and $Z(G)=1$, that is, both the Frattini subgroup of $G$ and the
center of $G$ are equal to 1.

Let $R$ be any non-nilpotent maximal subgroup of $G$ of order divisible by $p$. By the hypothesis,
$R\unlhd G$. Assume $N\nleq R$. Then $G=NR$. Note that $N\cap R\unlhd N$ and $N\cap R\unlhd R$,
one has $N\cap R\unlhd NR=G$. It follows that $N\cap R=1$ since $N\cap R<N$. Then $G=N\times R$,
which implies that $N\cong G/R$ is a cyclic group of prime order. This contradicts that $N$ is
non-solvable. Thus $N\leq R$.

Claim that $|N|$ is divisible by $p$. Otherwise, assume that $N$ has $p'$-order. Since $G$ is
non-solvable, one has that the intersection of all non-nilpotent maximal subgroups of $G$ is
equal to $\it\Phi$$(G)$ by {\rm\cite[Theorem 1]{szg}}. By above argument, the intersection of all
non-nilpotent maximal subgroups of $G$ is equal to 1. Since every non-nilpotent maximal subgroup
of $G$ of order divisible by $p$ contains $N$, there exists a non-nilpotent maximal subgroup $E$
of $G$ of $p'$-order such that $N\nleq E$. One has $G=NE$. It follows that $|G|=\frac{|N||E|}{|N\cap E|}$
has $p'$-order, a contradiction. Thus $|N|$ is divisible by $p$.

Since $Z(G)=1$, one has that all nilpotent maximal subgroups of $G$ are Sylow 2-subgroups of $G$
by $\rm\cite[Theorem 1]{rose}$. Then the maximal subgroups of $G$ may only be: Sylow 2-subgroups,
non-nilpotent maximal subgroups of order divisible by $p$ or non-nilpotent maximal subgroups of
$p'$-order. Assume that $q$ is a prime divisor of $|N|$ such that $q\neq p$. Let $P\in\Syl_p(N)$
and $Q\in\Syl_q(N)$. By Frattini argument, one has $G=NN_G(P)=NN_G(Q)$. Since $N$ is a minimal normal
subgroup of $G$ and $N$ is non-solvable, $N_G(P)<G$ and $N_G( Q)<G$.

It is obvious that both $N_G(P)$ and $N_G(Q)$ cannot be contained in any non-nilpotent maximal subgroup
of $G$ of order divisible by $p$ since every non-nilpotent maximal subgroup of $G$ of order divisible by $p$
contains $N$. (1) Assume $p=2$. Then $q$ is an odd prime. One has that $N_G(Q)$ can only
be contained in a non-nilpotent maximal subgroup $H$ of $G$ of $2'$-order. Then $G=NN_G(Q)=NH$. And
$N_G(P)$ can only be contained in a Sylow 2-subgroup $K$ of $G$. It follows that $G=NN_G(P)=NK$. One
has $|G|=|NH|=|NK|$. Then $\frac{|H|}{|H\cap N|}=\frac{|K|}{|K\cap N|}$. Note that $\frac{|H|}{|H\cap N|}>1$
is a $2'$-number, but $\frac{|K|}{|K\cap N|}>1$ is a $2$-power, a contradiction. (2) Assume $p>2$. One has
that $N_G(P)$ cannot be contained in any Sylow 2-subgroup of $G$ and $N_G(P)$ cannot be contained
in any non-nilpotent maximal subgroup of $G$ of $p'$-order, either, a contradiction.

Hence the counterexample of minimal order does not exist, then $G$ is solvable. \hfill\qed

\section{Proof of Theorem 1.4}\label{th4}

\medskip

\demo First we show that $G$ has a normal Sylow subgroup. Let $G$
be a counterexample of minimal order.

Claim $\it\Phi$$(G)=1$. Otherwise, assume $\it\Phi$$(G)\neq 1$. Since $|G/\it\Phi$$(G)|$ and $|G|$
have the same prime divisors, $G/\it\Phi$$(G)$ is also a group of order divisible by $p$ in which every
maximal subgroup is nilpotent or normal or has $p'$-order. By the minimality of $G$, one has that
$G/\it\Phi$$(G)$ has a normal Sylow subgroup $P\it\Phi$$(G)/\it\Phi$$(G)$, where $P$ is a Sylow
subgroup of $G$. Using Frattini argument, one gets that $G$ has a normal Sylow subgroup $P$,
a contradiction. Thus $\it\Phi$$(G)=1$.

Note that $G$ is solvable by Theorem 1.3. Let $N$ be a minimal normal subgroup of $G$, then $N$
is an elementary abelian group of prime power order. One has that there exists a maximal subgroup
$M$ of $G$ such that $N\nleq M$ since $\it\Phi$$(G)=1$. It follows that $G=NM$.

First consider the case when $M$ is a non-nilpotent maximal subgroup of $G$ of order divisible by
$p$. By the hypothesis, one has $M\unlhd G$. Note that $N\cap M$ is normal in $M$ and $N\cap M$
is also normal in $N$ as $N$ is abelian. Then $N\cap M$ is normal in $NM=G$. It follows that $N\cap M=1$
since $N$ is a minimal normal subgroup of $G$ and $N\cap M<N$. One has $G=N\times M$, which implies
that $M\cong G/N$. (1) Suppose that $G/N$ has $p'$-order. Then $N$ is a normal Sylow $p$-subgroup
of $G$, a contradiction. (2) Suppose that $G/N$ has order divisible by $p$. One has that $G/N$ has
a normal Sylow subgroup since the hypothesis of the theorem holds for $G/N$ and $|G/N|<|G|$. It
follows that $M$ has a normal Sylow subgroup $Q$ since $M\cong G/N$. $(i)$ Assume $(|N|,|Q|)\neq 1$.
Then $N\times Q$ is a normal Sylow subgroup of $G$, a contradiction. $(ii)$ Assume $(|N|,|Q|)=1$.
Then $Q$ is a normal Sylow subgroup of $G$, a contradiction, too.

Next consider the case when $M$ is a non-nilpotent maximal subgroup of $G$ of $p'$-order or a
nilpotent maximal subgroup of $G$ of $p'$-order. Note that $G=NM$ is a group of order divisible
by $p$ and $N$ is an elementary abelian group of prime power order. It follows that $N$ is a Sylow
$p$-subgroup of $G$, which is also a normal Sylow subgroup of $G$, a contradiction.

Finally consider the case when $M$ is a nilpotent maximal subgroup of $G$ of order divisible by
$p$. Assume $|N|=q^{m}$ for some prime $q$ and some positive integer $m\geq1$. Let $Q\in\Syl_q(M)$.
Then $NP\in\Syl_q(G)$. One has $NP\unlhd NM=G$, a contradiction.

By above arguments, the counterexample of minimal order does not exist and so $G$ has a normal
Sylow subgroup.

In the following we prove that $G$ has a Sylow tower.

Let $P_1$ be a normal Sylow $p_1$-subgroup of $G$.  (1) Suppose $p_1\neq p$. Observe that $G/P_1$
is also a group of order divisible by $p$ in which every maximal subgroup is nilpotent or normal or
has $p'$-order, arguing as above, one gets that $G/P_1$ has a normal Sylow $p_2$-subgroup $P_1P_2/P_1$,
where $P_2\in\Syl_{p_2}(G)$. (2) Suppose $p_1=p$. Let $M/P_1$ be a non-nilpotent maximal subgroup of
$G/P_1$, where $M$ is a non-nilpotent maximal subgroup of $G$ of order divisible by $p$. By the
hypothesis, one has $M\unlhd G$, then $G/P_1$ is a group in which every maximal subgroup is nilpotent
or normal. One can also get that $G/P_1$ has a normal Sylow $p_2$-subgroup $P_1P_2/P_1$ by
$\rm\cite[Theorem\,\,5]{shi}$, where $P_2\in\Syl_{p_2}(G)$.

Similarly, considering the quotient group $G/{P_1P_2}$, arguing as above we can get that $G/{P_1P_2}$
has a normal Sylow $p_3$-subgroup $P_1P_2P_3/{P_1P_2}$, where $P_3\in\Syl_{p_3}(G)$. And so on,
we can obtain a normal subgroups series:
\begin{equation}
P_1\unlhd{P_1P_2}\unlhd{P_1P_2P_3}\unlhd\cdots\unlhd{P_1P_2\cdots P_s}=G
\end{equation}
where $P_i\in\Syl_{p_i}(G)$ for $1\leq i\leq s$, which implies that $G$ has a Sylow tower. \hfill\qed

\section{Proof of Theorem 1.7}\label{th7}

\medskip

\demo Assume that $G$ is nilpotent, then $G$ is $q$-nilpotent and
$q$-closed for each prime divisor $q$ of $|G|$. In the following we assume that $G$ is non-nilpotent.
Since $G$ has Sylow tower by Theorem 1.4, $G$ can be written as $G=(P_1\times P_2\times\cdots\times P_s)
\rtimes (Q_1Q_2\cdots Q_t)$, where $P_i\in\Syl_{p_i}(G)$ and $P_i\unlhd G$ for $1\leq i\leq s$,
$Q_j\in\Syl_{q_j}(G)$ and $Q_j$ is not normal in $G$ for $1\leq j\leq t$. Let $M=P_1\times P_2\times
\cdots\times P_s$ and $N=Q_1Q_2\cdots Q_t$. Then $G=M\rtimes N$.

Considering the case when the Sylow $p$-subgroup of $G$ is normal in $G$. Assume $P_1\in\Syl_p(G)$, that is,
$p_1=p$.

(1) Suppose that $N$ is nilpotent. Then $G$ is obviously $p_i$-closed for every $1\leq i\leq s$ and $q_j$-nilpotent
for every $1\leq j\leq t$.

(2) Suppose that $N$ is non-nilpotent. Then there exists a maximal subgroup $N_0$ of $N$ such that $N_0$
is not normal in $N$. For the maximal subgroup $MN_0$ of $G$ of order divisible by $p$, $MN_0$ is not normal
in $G$. By the hypothesis, one has that $MN_0$ is nilpotent. Note that $|N:N_0|$ is a prime power by the
solvability of $G$. We can assume $Q_1\nleq N_0$ and $Q_j\leq N_0$ for every $2\leq j\leq t$. Then
$MN_0=P_1\times P_2\times\cdots\times P_s\times Q_{1_0}\times Q_2\times\cdots\times Q_t$, where
$Q_{1_0}<Q_1$. For any non-nilpotent maximal subgroup $H$ of $N$, one has that $MH$ is a non-nilpotent
maximal subgroup of $G$ of order divisible by $p$. By the hypothesis, $MH\unlhd G$. It follows that
$H\unlhd N$. That is every maximal subgroup of $N$ is nilpotent or normal, one has that $N$ has normal
Sylow subgroups by $\rm\cite[Theorem\,\,5]{shi}$. Assume $Q_j\unlhd N$ for some $2\leq j\leq t$, then
$N\leq N_G(Q_j)$. It follows that $G=N_G(Q_j)$ since $M\leq N_G(Q_j)$, this contradicts that $Q_j$ is
not normal in $G$ for every $2\leq j\leq t$. Thus $Q_j$ is not normal in $N$ for every $2\leq j\leq t$.
One has $Q_1\unlhd N$. Then $G=(P_1\times P_2 \times\cdots\times P_s)\rtimes (Q_1\rtimes (Q_2\times Q_3
\times\cdots\times Q_t))$. It is easy to see that $G$ is $p_i$-closed for every $1\leq i\leq s$, and
$G$ is $q_j$-nilpotent for every $2\leq j\leq t$. For $q_1$, $G$ is neither $q_1$-nilpotent
nor $q_1$-closed.

Next consider the case when the Sylow $p$-subgroup of $G$ is not normal in $G$. Assume $Q_1\in\Syl_p(G)$,
that is, $q_1=p$.

Here $G/M\cong N$ is a group of order divisible by $p$. Claim that $G/M$ is nilpotent. Otherwise, assume
that $G/M$ is non-nilpotent. Since $G/M$ is also a group of order divisible by $p$ in which every maximal
subgroup is nilpotent or normal or has $p'$-order, $G/M$ has a normal Sylow subgroup by Theorem 1.4,
which implies that $N$ has a normal Sylow subgroup. Assume $Q_j\unlhd N$ for some $1\leq j\leq t$.
Then $N\leq N_G(Q_j)$. It follows that $N_G(Q_j)$ has order divisible by $p$. Since $Q_j$ is not normal in
$G$ and all non-nilpotent maximal subgroups of $G$ of order divisible by $p$ are normal, $N_G(Q_j)$ can only
be contained in some nilpotent maximal subgroup of $G$ of order divisible by $p$. It implies that $N$ is
nilpotent, this contradicts that $N\cong G/M$ is non-nilpotent. Hence $N\cong G/M$ is nilpotent. Then
$G$ is $p_i$-closed for every $1\leq i\leq s$ and $q_j$-nilpotent for every $1\leq j\leq t$.\hfill\qed

\section{Proofs of Theorems 1.8, 1.9 and 1.10}\label{th8}

\medskip

\f {\bf Proof of Theorem 1.8.}\ \ \ Let $G$ be a group in which every maximal subgroup is nilpotent
or normal. Suppose that $G$ is nilpotent, then $G$ is $q$-nilpotent and $q$-closed for each prime
divisor $q$ of $|G|$.

Next we suppose that $G$ is non-nilpotent. Note that $G$ has Sylow tower by {\rm\cite[Theorem 5]{shi}},
one has $G=(P_1\times P_2\times\cdots\times P_s)\rtimes (Q_1Q_2\cdots Q_t)$, where
$P_i\in\Syl_{p_i}(G)$ and $P_i\unlhd G$ for $1\leq i\leq s$, $Q_j\in\Syl_{q_j}(G)$ and $Q_j$ is not
normal in $G$ for $1\leq j\leq t$. Let $K=P_1\times P_2\times\cdots\times P_s$ and $L=Q_1Q_2\cdots Q_t$.
Then $G=K\rtimes L$.

Claim that $L$ is nilpotent. Otherwise, assume that $L$ is non-nilpotent. For every non-nilpotent maximal
subgroup $L_0$ of $L$, one has that $KL_0$ is a non-nilpotent maximal subgroup of $G$. By the hypothesis,
one has $KL_0\unlhd G$. It follows that $L_0\unlhd L$. That is, $L$ is a group in which every maximal
subgroup is nilpotent or normal. By {\rm\cite[Theorem 5]{shi}}, $L$ has a normal Sylow subgroup. Assume
$Q_1\unlhd L$. Then $L\leq N_G(Q_1)<G$. Since $L$ is non-nilpotent, $N_G(Q_1)$ can only be contained in
a non-nilpotent maximal subgroup $H$ of $G$. By the hypothesis, $H\unlhd G$. Then by Frattini argument,
one has $Q_1\unlhd G$, a contradiction. Thus $L$ is nilpotent. It follows that $G$ is $p_i$-closed
for every $1\leq i\leq s$ and $q_j$-nilpotent for every $1\leq j\leq t$.\hfill\qed

\medskip

\f {\bf Proof of Theorem 1.9.}\ \ \ If $G$ is nilpotent, it is obvious that $G$ is $q$-nilpotent and
$q$-closed for each prime divisor $q$ of $|G|$. In the following assume that $G$ is non-nilpotent.
Arguing as in proof of Theorem 1.7, one has $G=(P_1\times P_2\times\cdots\times P_s)
\rtimes (Q_1Q_2\cdots Q_t)$, where $P_i\in\Syl_{p_i}(G)$ and $P_i\unlhd G$ for $1\leq i\leq s$,
$Q_j\in\Syl_{q_j}(G)$ and $Q_j$ is not normal in $G$ for $1\leq j\leq t$. Assume $M=P_1\times
P_2\times\cdots\times P_s$ and $N=Q_1Q_2\cdots Q_t$. That is, $G=M\rtimes N$.

We will show that $N$ is nilpotent.

First assume $p_1=p$, that is, the Sylow $p$-subgroup of $G$ is normal and $P_1\in\Syl_p(G)$.
Let $N_0$ be any maximal subgroup of $N$. Then $M\rtimes N_0$ is a maximal subgroup of $G$ of
order divisible by $p$. By the hypothesis, $M\rtimes N_0$ is nilpotent. It follows that $N_0$
is nilpotent. Then $N$ is nilpotent or minimal non-nilpotent. Claim that $N$ cannot be minimal
non-nilpotent. Otherwise, assume that $N$ is a minimal non-nilpotent group. We can assume
$N=Q_1\rtimes Q_2$ by {\rm\cite[Theorem 9.1.9]{rob}}, where $Q_1\unlhd N$. Let $Q_2'$ be
any maximal subgroup of $Q_2$. Then $M\rtimes(Q_1\rtimes Q_2')$ is a maximal subgroup of $G$
of order divisible by $p$. By the hypothesis, $M\rtimes(Q_1\rtimes Q_2')$ is nilpotent. It
follows that $M\leq N_G(Q_1)$. Note that $N\leq N_G(Q_1)$. Then $N_G(Q_1)=MN=G$, which
implies that $Q_1\unlhd G$, a contradiction. Thus $N$ is nilpotent.

Second assume $q_1=p$, that is, the Sylow $p$-subgroup of $G$ is not normal and $Q_1\in\Syl_p(G)$.
It is obvious that $N<G$ and $N$ has order divisible by $p$. Then there exists a maximal
subgroup $L$ of $G$ of order divisible by $p$ such that $N\leq L$. By the hypothesis, $L$
is nilpotent. It follows that $N$ is also nilpotent.

By the nilpotence of $N$, one has that $G$ is $p_i$-closed for every $1\leq i\leq s$ and
$q_j$-nilpotent for every $1\leq j\leq t$.
\hfill\qed

\medskip

\f {\bf Proof of Theorem 1.10.}\ \ \ It is clear that the result holds if $G$ is nilpotent.
Considering the case that $G$ is non-nilpotent. Arguing as above, let $G=(P_1\times P_2
\times\cdots\times P_s)\rtimes (Q_1Q_2\cdots Q_t)$, where $P_i\in\Syl_{p_i}(G)$ and
$P_i\unlhd G$ for $1\leq i\leq s$, $Q_j\in\Syl_{q_j}(G)$ and $Q_j$ is not normal
in $G$ for $1\leq j\leq t$. Let $M=P_1\times P_2\times\cdots\times P_s$ and $N=Q_1Q_2
\cdots Q_t$. One has $G=M\rtimes N$.

Claim that the Sylow $p$-subgroup of $G$ is normal. Otherwise, assume that the Sylow
$p$-subgroup of $G$ is not notmal. Let $Q_1\in\Syl_p(G)$. Note that $N_G(Q_1)<G$ and
$N_G(Q_1)$ has order divisible by $p$. Then there exists a maximal subgroup $R$ of
$G$ of order divisible by $p$ such that $N_G(Q_1)\leq R$. By the hypothesis, $R\unlhd G$.
It follows that $Q_1\unlhd G$ by Frattini argument, a contradiction.

Thus the Sylow $p$-subgroup of $G$ is normal. Let $P_1\in\Syl_p(G)$. For every maximal
subgroup $N_0$ of $N$, $MN_0$ is a maximal subgroup of $G$ of order divisible by $p$.
By the hypothesis, $MN_0\unlhd G$. It follows that $N_0\unlhd N$ and then $N$ is nilpotent.
Thus $G$ is $p_i$-closed for every $1\leq i\leq s$ and $q_j$-nilpotent for every
$1\leq j\leq t$.\hfill\qed

\bigskip

\f {\bf Acknowledgements}

\medskip

The authors are thankful to everyone who provides valuable suggestions and helpful comments
for improving our paper.

\end{document}